\documentclass[conference]{IEEEtran}

\usepackage{latexsym,amssymb,amsmath,theorem,stmaryrd}

\usepackage{balance}
\usepackage[draft,inline,nomargin]{fixme}
\usepackage{graphicx}
 \usepackage{times}
\usepackage[usenames]{color}
 \usepackage{multirow}

\newtheorem{example}{Example}

\newcommand{\C}{\mathcal C}

\newcommand{\Y}{\mathbf{W}}
\newcommand{\Q}{\mathbf{Q}}
\newcommand{\Qb}{\overline\Q}
\newcommand{\E}{\mathbf{E}}
\newcommand{\I}{\mathbf{I}}
 \newcommand{\win}{\C}
\newcommand{\nwin}{\bar \C}
\newcommand{\s}{{\bf{\pi}}}

\title{Bounding the Equilibrium Distribution of\\  
Markov Population Models}

\author{\IEEEauthorblockN{Tu\v{g}rul Dayar$^1$ \qquad \qquad Holger
    Hermanns  \qquad \qquad 
David Spieler \qquad  \qquad Verena Wolf\\
Saarland University -- Computer Science, Saarbr\"{u}cken, Germany
}}
\IEEEspecialpapernotice{(Extended Abstract)}

\begin{document}
 \maketitle

\renewcommand{\thefootnote}{\fnsymbol{footnote}}
\footnotetext[1]{%
  T. Dayar is
currently on sabbatical
 leave from Bilkent University, Turkey.}
\renewcommand{\thefootnote}{\arabic{footnote}}

\section{Introduction}
 
Markov   population models (MPMs) are  continuous-time 
Markov chains (CTMCs) that describe the dynamics and interactions
of different populations.
They have  important applications in the life science domain, in particular 
  in ecology, epidemics, and biochemistry.
Depending on the system under study, a member of a population 
represents an individual that belongs to a certain  biological species,
an organism that suffers from an infectious disease, or
a certain type of molecule in a living cell.
Thus, if $n$ is the number of populations, the state space of the MPM 
is  $\mathbb N^n$, that is, a state is a vector $(x_1,\ldots,x_n)$ 
of non-negative integers, where the entry $x_i$ is the size of the $i$-th 
population. 
Typically, the transitions of an MPM are described by a finite set of
transition classes such that each class specifies a (possibly infinite) 
number of edges in the underlying transition graph. 
For instance, we may have one class to represent the death of individuals. 
 In biochemistry, each chemical reaction describes a class of transitions 
in the associated MPM.
Often, the corresponding transition rates are state-dependent, e.g., 
the rate at which individuals of a certain population die may depend 
on the population size.

The structural regularity of MPMs often enables accurate
approximations of the system behavior. One such example is the
widely-used deterministic approximation of the dynamics of chemical
reaction networks~\cite{kurtz} that represents the states as a
continuum. But if one or more populations are small a discrete
representation of the population sizes is important and continuous
approximations are inaccurate. This effect has also been observed
experimentally in the context of chemical
reactions~\cite{Thattai,swain,paulsson}.  In such cases the analysis
of the MPM becomes difficult. Closed-form solutions are only possible
in special circumstances~\cite{Jahnke} and numerical solution
techniques suffer from the problem that a very large or even infinite
state space has to be explored. Therefore, Monte-Carlo simulation is
in widespread use to estimate transient or stationary measures of the
MPM.

Recently, progress has been made on numerically approximating the
transient distribution of an MPM at particular time
instances~\cite{HIBI09,sparsegrid,Inexact}.  These approaches exploit
the fact that only a subset of the state space is needed to give an
accurate approximation and that only a small amount of probability
mass is located above a certain population threshold. The intuitive
explanation for this is simple since within a fixed time interval it
is extremely unlikely that the populations reach certain thresholds.
For the long-run behavior of the system, however, this argument does
not hold since it is a priori not known if and where the system
stabilizes.

\begin{figure}[b!]
   \includegraphics[width=7.5cm]{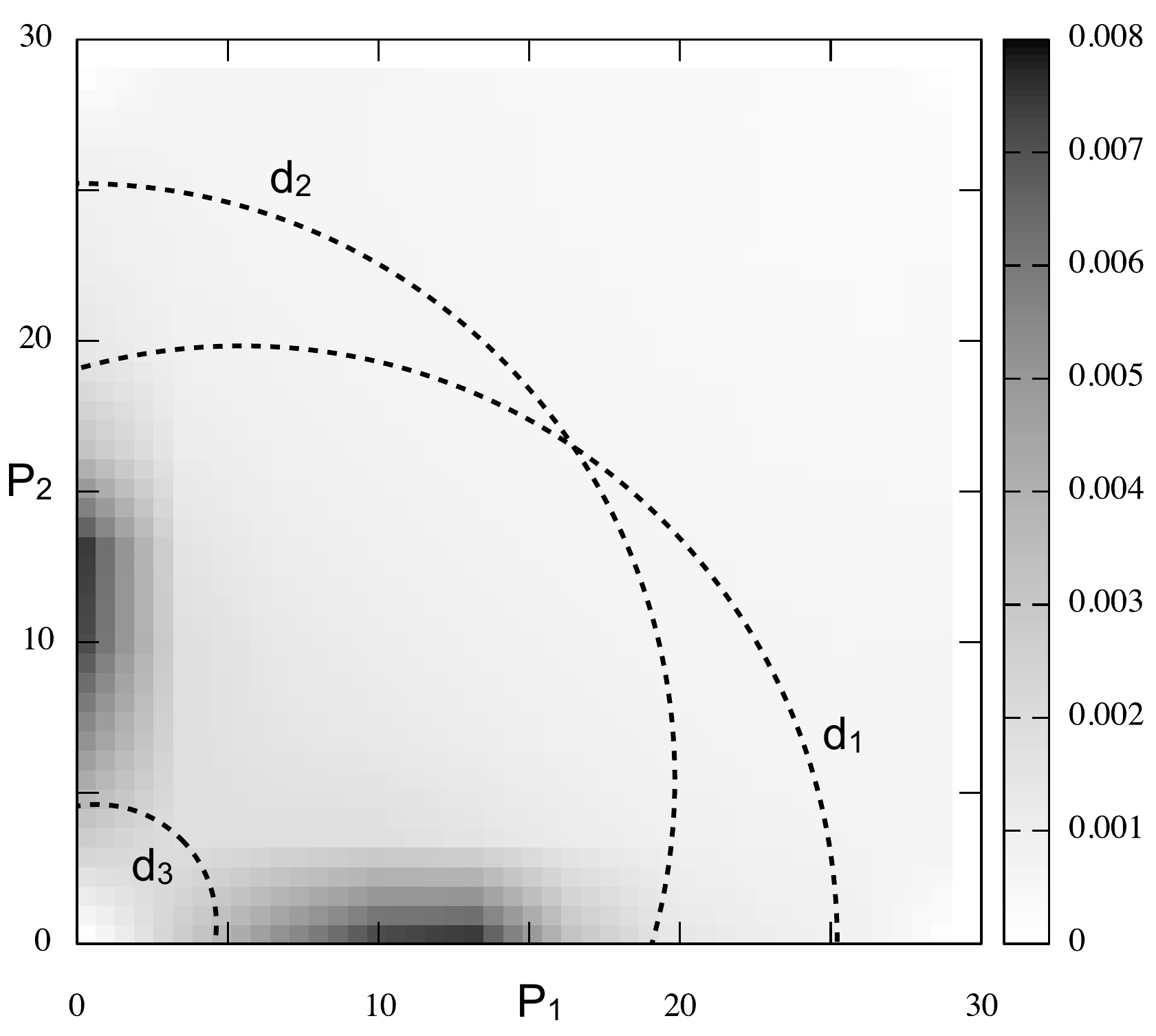} 
                    \caption{The equilibrium distribution $\pi$ of an MPM 
that describes a gene regulatory network and geometric bounds enclosing
a subset $\C$ of the state space with $\sum_{c\in\C}\pi(c) > 0.9$. }%
        \label{fig:exswitch}%
   \end{figure}

In this paper, we consider the problem of computing an accurate 
approximation of the equilibrium distribution of an MPM.
Assuming that the MPM is ergodic, we first derive geometric bounds
for the equilibrium distribution, i.e., we find those regions of the 
state space, where most of the probability mass is located in the limit.
Then we perform a local refinement for these regions in order to 
bound the probabilities of individual states.

Let $\{X(t),t\ge 0\}$ be an MPM.  The first step relies on an analysis
of a \emph{drift function} $d$ that associates with each state $x$ the
expected change $d(x)=\frac{d}{dt}E[g(X(t))\mid X(t)=x]$.  Here, $g$ is a
function that associates with each state $x$ a non-negative real
number, also called Lyapunov function.  We illustrate this by means of an example.
Figure~\ref{fig:exswitch} shows the plot of the equilibrium
distribution of an MPM that describes a gene regulatory network.  The
system is bistable, that is, in the long-run the probability mass is
concentrated at two distinct regions in the state space.  In these
regions the drift $d(x)$ is maximal. We determine geometric bounds
(here depicted as dashed lines) by using a simple threshold on the
drift, i.e., we consider the set $\C$ of all states where the drift is greater 
than a certain threshold. Then $\C$ contains those states, where most of the probability
mass is located in equilibrium.  For the amount of probability within
$\C$, we derive tight bounds.

In the second step of our approach, we consider 
the stochastic complement of $\C$~\cite{meyer} which allows us to derive 
bounds for the probabilities of individual states. 
The stochastic complement is a finite CTMC with state space $\C$,
where each outgoing transition leading to a state not in $\C$ is redirected
to $\C$. This redirection is defined in such a way that the equilibrium distribution
of the stochastic complement gives the conditional equilibrium probabilities 
for the states of $\C$. Since the exact redirection probability can only
 be obtained from a full solution of the infinite system, we consider 
candidates that give upper and lower bounds. 
Together with the first step of the approach, this yields bounds for the 
equilibrium probabilities of all states of the MPM $X$.

\emph{The full paper version will contain all relevant proof details.} 
 
\section{Markov Population Models}
We consider  a class of time-homogeneous CTMCs that can be described 
by a finite set of transition classes. A transition class $\tau$ 
is a pair $(\alpha,v)$, where $\alpha:\mathbb N^n\to \mathbb R_{\ge 0}$
is a function that determines the transition rate and 
  $v\in \mathbb Z^n$ is a change vector that determines the successor state 
of the transition. Thus, if $x\in\mathbb N^n$ and
$\alpha(x)>0$ then there is a transition  from state $x$ 
to state $x+v$ with rate $\alpha(x)$. 
We assume that $v$ has at least one nonzero entry and that 
$\alpha(x)$ is a polynomial in $x=(x_1,\ldots,x_n)$. 
Let $\{\tau_1,\ldots,\tau_k\}$ be a set of transition classes
with distinct change vectors.
Let $\Q$ be the infinite matrix such that the entry $\Q(x,x+v_j)$ equals 
$\alpha_j(x)$, where $\tau_j=(\alpha_j,v_j)$ and $j\in\{1,\ldots,k\}$.
If we define the diagonal entries of $\Q$ as the negative sum of the off-diagonal 
entries, then $\Q$ is the infinitesimal generator of a 
  CTMC $\{X(t),t\ge 0\}$. 
The matrix $\Q$ has only finitely many nonzero entries in each row and in each column, i.e.,
$\Q$ is an infinite matrix with finite rows/columns.
Note that $\sup_{x}|\Q(x,x)|$ may   be infinite and that
  the number of  states reachable from a given initial state may be infinite.

\section{Geometric Bounds}
\label{sec:geometric_bounds}
We assume that $\tau_1,\ldots,\tau_k$ are such that $X$ is 
ergodic and $\pi$ is the equilibrium distribution of $X$.
 Then there exists a  
function $g^*:\mathbb N^n\to\mathbb R_+$, a finite subset $\C$, and a constant $\gamma>0$ such that~\cite{tweedie}
\begin{equation}
 \label{eq:conditions}
   \begin{array}{l@{\,}c@{\ }l}
(i)&  \frac{d}{dt}E[g^*(X(t))\mid X(t)=x]\le -\gamma, &\forall x\in \mathbb N^n\setminus \C,\\[1ex]
 (ii)&   \frac{d}{dt}E[g^*(X(t))\mid X(t)=x]<\infty, &\forall x\in   \C,\\[1ex]
(iii) & \multicolumn{2}{l}{ \{x\in \mathbb N^n\mid g^*(x)\le \ell\} \mbox{ is finite for any } \ell<\infty.}
   \end{array}
\end{equation}
Let $c$ be a positive number with $c\ge \max_{x\in \C}E[g^*(X(t))\mid X(t)=x]$ and 
let $g(x)=\frac{g^*(x)}{\gamma+c}$. In the sequel, we will refer to $g$ as the Lyapunov
function.
The first two conditions in Eq.~\eqref{eq:conditions} are now equivalent to
\begin{equation}\label{eq:glynn}
\begin{array}{lcl}
 \frac{d}{dt}E[g(X(t))\mid X(t)=x]
&\le& \frac{c}{c+\gamma} -   \bar \chi_{\C}(x),
\end{array}
\end{equation}
where   $\bar\chi_{\C}(x) =1$
if $x\not\in \C$ and 0 otherwise.
Note that if we multiply Eq.~\eqref{eq:glynn}
 with $\pi(x)$ and sum  over $x$, the left-hand side 
becomes zero and we arrive at 
$$
\sum_{x\not\in \C}\pi(x)\le \frac{c}{c+\gamma}.
$$
Thus, we can use Eq.~\eqref{eq:glynn} to bound the probability mass 
outside of $\C$. 

For a given Lyapunov function $g$, we define  the drift $d(x)$ in state $x$ as
$d(x)=\frac{d}{dt}E[g(X(t))| X(t)=x]$. 
Since $X$ has a transition class description,  
$$d(x)= \sum_{j=1}^k \alpha_j(x)(g(x+v_j)-g(x)).$$
 If $g$ is a polynomial of degree $\ell_g>0$  
and the highest degree of the rate functions $\alpha_j$ is $\ell_d$  
then the drift function $d$ 
is at most a polynomial of degree $(\ell_g-1) \ell_d$.
For most population models, $\ell_d\le 2$. Moreover, often a degree of $\ell_g=2$
is sufficient for $g$. In these cases, we can easily determine the maxima of $d$
and use Eq.~\eqref{eq:glynn} to derive bounds for the probability mass outside of  $\C$.
 Note that we can also bound the probability mass inside $\C$ with symmetric arguments~\cite{GlynnZeevi}.

\begin{example}
\label{ex:exswitch}
 We consider a gene regulatory
network called the exclusive switch~\cite{exswitch}.
It consists of two genes with promotor regions.
Each of the two gene products $P_1$ and $P_2$ inhibits the 
expression of the other product if a molecule is bound to the 
respective promotor region. More precisely, 
if the promotor region is free, molecules of  both types 
$P_1$ and $P_2$ are produced.
If a molecule of type $P_1$ is bound to the promotor 
region of $P_2$, only molecules of type $P_1$ are produced.
If a molecule of type $P_2$ is bound to the promotor 
region of $P_1$, only molecules of type $P_2$ are produced.
No other configuration of the promotor region exists.
The system has six chemical species of which two 
have an infinite range, namely $P_1$ and 
$P_2$. We define the  transition classes 
$\tau_j=\left(\alpha_j,v_j\right)$, $j\in\{1,\ldots,8\}$ as follows.
 	\begin{itemize}
		 \item For $j\in\{1,2\}$ we describe production of  $P_j$
by  $v_j(x)=e_j$ and 
$\alpha_j(x)=0.05  x_{2+j}$. Here, $x_{2+j}$ the number of active genes that produce $P_j$,
which is either zero or one.
The $j$-th entry of a state $x$ represents the number of $P_j$ molecules and
the vector  $e_j$
is such that all its entries are zero except the $j$-th entry which is one.

\item We describe degradation of  $P_j$
by $v_{j+2}=-e_j$ and 
$\alpha_{j+2}(x)=0.005  x_j$.  Here, $x_j$ denotes the number of $P_j$ molecules.

\item We model the binding of $P_1$  to the promotor (which inhibits the gene 
that is responsible for the production of $P_2$) as
  $v_{5}=-e_1-e_{4}+e_{6}$, and 
$\alpha_{5}(x)=0.1  x_1 x_3  x_4$. 
Here, $x_{6}$ is one if a molecule of type $P_{1}$ 
is bound to the promotor region   and zero otherwise.
Note that $\alpha_5$ is zero in all states where the promotor is not free ($x_3=0$
or $x_4=0$).
\item We model the binding of $P_2$  to the promotor (which inhibits the gene 
that is responsible for the production of $P_1$) as
  $v_{6}=-e_2-e_{3}+e_{5}$, and 
$\alpha_{6}(x)=0.1  x_2 x_3  x_4$. 
Here, $x_{5}$ is one if a molecule of type $P_{2}$ 
is bound to the promotor region   and zero otherwise.
\item For unbinding of $P_1$ we define
  $v_{7}=e_{1}+e_{4}-e_{6}$, and 
$\alpha_{7}(x)=0.005 x_{6}$.
\item For unbinding of $P_2$ we define
  $v_{8}=e_{2}+e_{3}-e_{5}$, and 
$\alpha_{8}(x)=0.005 x_{5}$.
 \end{itemize}
We  use the Lyapunov  function $g$ given by
$g(x) = x_1^2+x_2^2+\ldots+ x_6^2.$
Consequently, the drift becomes
$$d(x) = 0.05  x_3  (2 x_1 + 1) + 0.05  x_4  (2 x_2 + 1)$$
$$+~0.1  x_3  x_4  x_1  ( -2  x_4 -2  x_1 + 2  x_5 + 3)$$
$$+~0.1  x_3  x_4  x_2  ( -2  x_3 -2  x_2 + 2  x_6 + 3)$$
$$+~0.005  x_1  (-2  x_1 + 1) + 0.005  x_2  (-2  x_2 + 1)$$
$$+~0.005  x_5  (-2  x_5 + 2  x_3 + 2  x_2 + 3)$$
$$+~0.005  x_6  (-2  x_6 + 2  x_4 + 2  x_1 + 3).$$
With the initial condition $x_{j+2} = 1$, invariantly it holds that
$x_{j+2} \in \{0,1\}$ and $x_{j+2} = 1 - x_{j+4}$. The global maxima of $d(x)$ 
(when considering real valued $x_{j}$) therefore are found
at $x^{(m1)} = (0.25,5.75,0,1,1,0)$ and
$x^{(m2)} = (5.25,0.25,0,1,1,0)$. 
The maximal value of the drift   in the
reachable part of the state space consequently is lower or equal to
$$c = d(x^{(m1)}) = d(x^{(m2)}) = 0.38625.$$

We are interested in a set $\C$ with $\sum_{x\in\C}\pi(x) > 1-\epsilon$,
where $\epsilon\in (0,1)$ is an a priori chosen threshold. 
Let $\gamma>0$ be such that $\epsilon=c/(c+\gamma)$.
We choose $\C$ such that is contains all states where 
the drift is greater than $\gamma$. 
Note that $\C$ is finite since $g$ fulfills condition $iii)$ in Eq.~\eqref{eq:conditions}.
It is easy to verify that 
Eq.~\eqref{eq:glynn} holds if 
we scale the Lyapunov function $g$ by $\gamma+c$, that is,
$g(x) =(x_1^2+x_2^2+\ldots+ x_6^2)/(\gamma+c).$
We retrieve 
the normalized drift
$$d_s(x) = \frac{1}{\gamma+c} d(x).$$
Therefore, $\C = \{x \in\mathbb{N}^6 ~|~ d_s(x) > \epsilon-1\}$ and
we get the desired bound for the equilibrium probability inside $\C$.
With the constraints on $x_j$ and $x_{j+2}$ we only have to consider the bounds for $x_1$
and $x_2$ and derive three cases, namely:
\begin{enumerate}
  \item $d_1(x) = d_s(x_1,x_2,1,0,0,1) = -0.01 x_1^2 - 0.1 x_2^2 + 0.115 x_1 + 0.005 x_2 + 0.055 > \epsilon - 1$,
  \item $d_2(x) = d_s(x_1.x_2,0,1,1,0) = -0.01 x_1^2 - 0.1 x_2^2 + 0.005 x_1 + 0.115 x_2 + 0.055 > \epsilon - 1$,
  \item $d_3(x) = d_s(x_1,x_2,1,1,0,0) = -0.21 (x_1^2 + x_2^2) + 0.205 (x_1+x_2) + 0.1 > \epsilon - 1$,
\end{enumerate}
illustrated by Figure \ref{fig:exswitch} for $\epsilon=0.1$.
\end{example}
 
In order to apply the approach described above, 
one has to find  a Lyapunov function $g$ such that Eq.~\eqref{eq:glynn}
holds for some $c$ and $\gamma$. 
This may become difficult for complex systems even though often a 
quadratic function is sufficient and it is possible to optimize  
the coefficients. 

\section{Conditional Probabilities}
\label{sec:cond_prob}
In order to derive probability bounds for the individual states in $\C$, we
consider the following 
 partitioning 
$$\Q = \begin{bmatrix} \Q[\win,\win] &\quad & \Q[\win,\nwin]  \\ \\ \Q[\nwin,\win] &\quad& \Q[\nwin,\nwin] \end{bmatrix}$$
of the generator matrix $\Q$
into blocks that describe the transitions within $\C$, to $\bar \C$, within $\bar \C$
and back to $C$.
Let $\E$ denote the embedded Markov chain of $\Q$, i.e.,
$$\E = \I - D_{\Q}^{-1} \Q = \begin{bmatrix} \E[\win,\win] &\quad & \E[\win,\nwin]  \\ \\ \E[\nwin,\win] &\quad& \E[\nwin,\nwin] \end{bmatrix}\!\!\!,$$
where $\I$ denotes the identity matrix and 
$D_{\Q}$
is a diagonal matrix such that   $D_{\Q}(x,x)=-\Q(x,x)$ for all states $x$.

The \emph{stochastic complement} of $\win$ is defined as
$$\Qb = \Q[\win,\win] + \Q[\win,\nwin] \sum_{i=0}^{\infty}(\E[\nwin,\nwin])^i \E[\nwin,\win].$$
Using the fact that $X$ is ergodic, we are able to show that
  $\Qb$  is   well-defined and  is
 the infinitesimal generator of a finite ergodic
  CTMC.

Let $\s_{\win}$ denote  the equilibrium distribution of  $\Qb$.
For finite discrete-time Markov chains, it has been shown in the seminal work by Meyer~\cite{meyer} 
that the entries of $\s_{\win}$ are equal to the 
conditional equilibrium probabilities of $X$, i.e., 
$\s_{\win}(x) =  {\s(x)}/{\sum_{c\in\C}\s(c)}$ for all $x$.
In what follows, we extend this result to infinite MPMs.

 The construction of the stochastic complement requires that transition probabilities
inside the  infinite set $\nwin$ are calculated. 
Since this is infeasible for MPMs, we apply a similar technique as proposed by
 Courtois and Semal  for finite CTMCs\cite{courtois,courtoissemal}.
The idea is to only consider the set $\C$ and redirect the transitions that lead 
from $\C$ to $\nwin$ back to $\C$. The matrix $\Qb$ contains the ``exact'' redirections
of the transitions, i.e., the solution of the corresponding CTMC gives the conditional 
probabilities of the states in $\C$. Since we cannot construct  $\Qb$, we redirect the 
transitions in such a way that we obtain upper and lower bounds.  

We first consider the substochastic matrix $\Y$ given by
$$\Y = \I + \frac{1}{\lambda}\Q[\win,\win]$$
with $\lambda > \max_{x\in \C} -\Q[\win,\win](x,x)$.
If we increase the $j$-th column of $\Y$ such that it becomes a
stochastic matrix, it is easy to see that the result represents an
ergodic discrete-time Markov chain for $j = 1,\ldots,|\C|$. 
When computing the conditional probabilities for relatively small 
values of $|\C|$, one can pass the slack probability mass summing up the 
rows of $\Y$ to 1 in an extra column to a dummy state and add an extra 
row corresponding to the dummy state which redirects the system to state 
$j$ with probability 1. After removing the redundant last equation, the 
transposed linear system can be written such that $\Y^T-I$ is the 
coefficient matrix and $e_j$ is the right-hand side vector. This makes 
it possible to LU factorize $\Y^T-I$ only once, and obtain the solution 
by forward and backward substitutions followed by normalization for 
$j = 1,\ldots,|\C|$. Now,
let $\s^{\Y}_j$ be the associated equilibrium distribution. From this, 
we are able to derive that for all $x \in \C$
$$\min_j \s^{\Y}_j(x) \leq \frac{\s(x)}{\sum_{c\in\win}\s(c)} \leq \max_j \s^{\Y}_j(x),$$

For a given threshold $\epsilon>0$, we first determine the set $\C$
as described in
Section \ref{sec:geometric_bounds}. Then
we bound the individual (unconditioned) state probability of a state $x\in\C$  by
$$(1-\epsilon)\min_j \s^{\Y}_j(x) \leq \s(x) \leq \max_j \s^{\Y}_j(x).$$
 
%
%
\begin{example}
For Example \ref{ex:exswitch}
we received tight  bounds
on the individual conditioned 
equilibrium probabilities inside $\C$ for $\epsilon = 0.1$  (cf. Fig. \ref{fig:exswitch}). 
In  Fig. \ref{fig:diff_bounds} we plot the difference between upper and lower bounds
for all states. We achived a precision of $\delta = 3.5\cdot10^{-4}$, i.e., the maximal difference between
upper and lower bound was $\delta$.
Note that a lower bound can be retrieved by multiplying the distribution
of Fig. \ref{fig:exswitch} by $1-\epsilon=0.9$. The upper bound of $\pi(x)$ is given by
$\max_j\pi_j^\Y(x)$.
There is a total of 1671 states in $\C$ for the chosen value of $\epsilon$.

\begin{figure}
\begin{center}
  \includegraphics[width=0.4\textwidth]{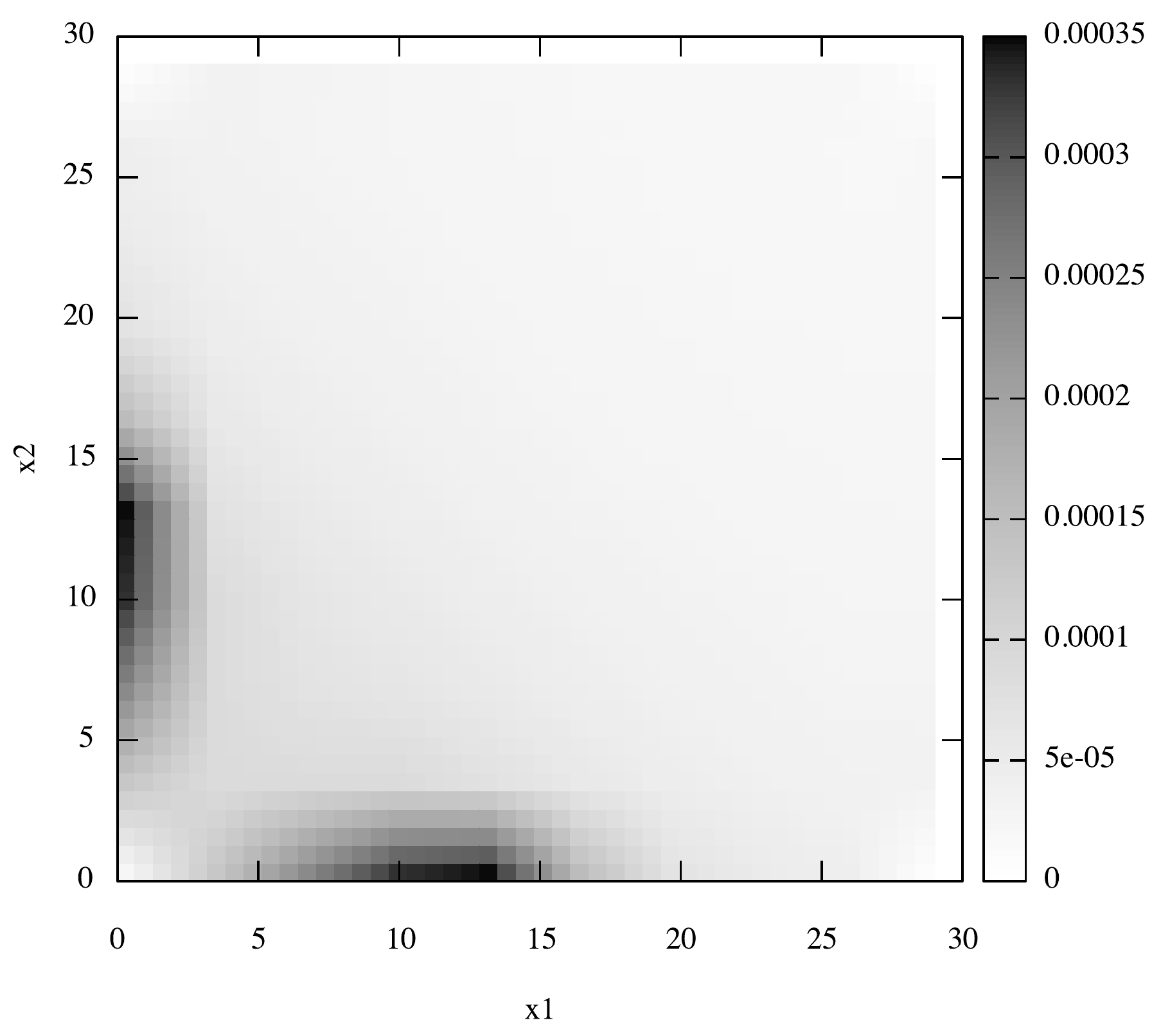}
  \caption{Difference between upper and lower bounds on states in $\C$ for $\epsilon=0.1$.}
    \label{fig:diff_bounds}
\end{center}
\end{figure}
\end{example}
\section{Conclusion}
We have demonstrated how to calculate equilibrium probability bounds of
infinite MPM by combining Lyapunov theory with numerical approximation
and bounding techniques. Much remains to be done with respect to
implementation efficiency, since various time-space tradeoffs appear
worthwhile to be explored. 

\newpage
  \bibliographystyle{abbrv}
 \bibliography{MPM}

\end{document}